\documentclass{article}
\usepackage[utf8]{inputenc}
\usepackage{graphicx}
\usepackage{authblk}
\usepackage{biblatex}
\usepackage{amsmath}
\usepackage{amsfonts}
\usepackage{color}
\usepackage{bm}
\usepackage{tabularx}
\usepackage{url}

\title{Homochiral inflation for the aperiodic monotile Tile(1,1)}
\author{Marianne Imperor-Clerc, Jean-François Sadoc}
\affil{Laboratoire de Physique des Solides, CNRS and Université Paris-Saclay, 91400 Orsay, France email:~ marianne.imperor@cnrs.fr}
\date{\today}
\begin{document}
\maketitle
\begin{abstract}
The recently discovered chiral monotile Tile(1,1) is tiling the plane in a quasiperiodic fashion by taking twelve different orientations when applying $2\pi/12$ rotation. An homochiral inflation construction
of such a quasiperiodic tiling is proposed where the chirality of the monotile is completely fixed at all inflation steps, avoiding to exchange its chirality between two successive steps. Doing so, the twelve possible orientations of the monotile are explicitly coded and the key difference between odd and even orientations is taken into account. The tiling is decomposed using only two different clusters, $\Gamma$ and $\Omega$, each of them taking six possible orientations. This gives a total set of twelve metatiles, which assembly can be mapped onto a triangular lattice. This approach allows to properly separate rotation and translation symmetry elements relating monotiles together. As all possible orientations of the two clusters are already incorporated in the twelve metatiles, positions of adjacent metatiles are given by translations which are along three equivalent directions ($2\pi/3$ rotation) as evidenced by junction lines. Finally, thanks to the homochiral inflation, the orientation distribution of the monotile at each inflation step is computed.   
\end{abstract}

\section{Introduction}

Since its discovery in May 2023 \cite{smith2023chiral}, the chiral monotile called Tile(1,1) or also Spectre has attracted a large interest in the scientific community and has reached rapidly the general public audience. This chiral monotile is solving the long standing 'Ein Stein' problem \cite{Baake_sym_2012} by giving rise to a quasiperiodic tiling of the plane by a single polygonal tile without allowing its mirror reflection, contrary to the related Hat monotile \cite{smith2023aperiodic}. In such tiling, the chiral monotile is taking 12 different orientations by applying $2\pi/12$ rotation. An inflation construction has been given, where the chirality of the monotile is changing at each inflation step \cite{smith2023chiral,  voss_2024_algorithm_mathlab}. However, in order to better understand the quasiperiodic properties of this tiling, there is a need to incorporate explicitly the orientation distribution of the monotile in this inflation construction and this is the main goal of this paper. We propose an homochiral inflation construction of a quasiperiodic tiling with Tile(1,1) where the chirality of this monotile is entirely fixed. Furthermore, rotation and translation symmetry elements are separately taken into account in this description. It allows to derive the orientation distribution of the monotile and the proportions of the 12 orientations of the monotile are computed at each inflation step.

\section{Tiles and clusters}
\label{sec:tiles_clusters}
\begin{figure*}
\includegraphics[width=\textwidth]{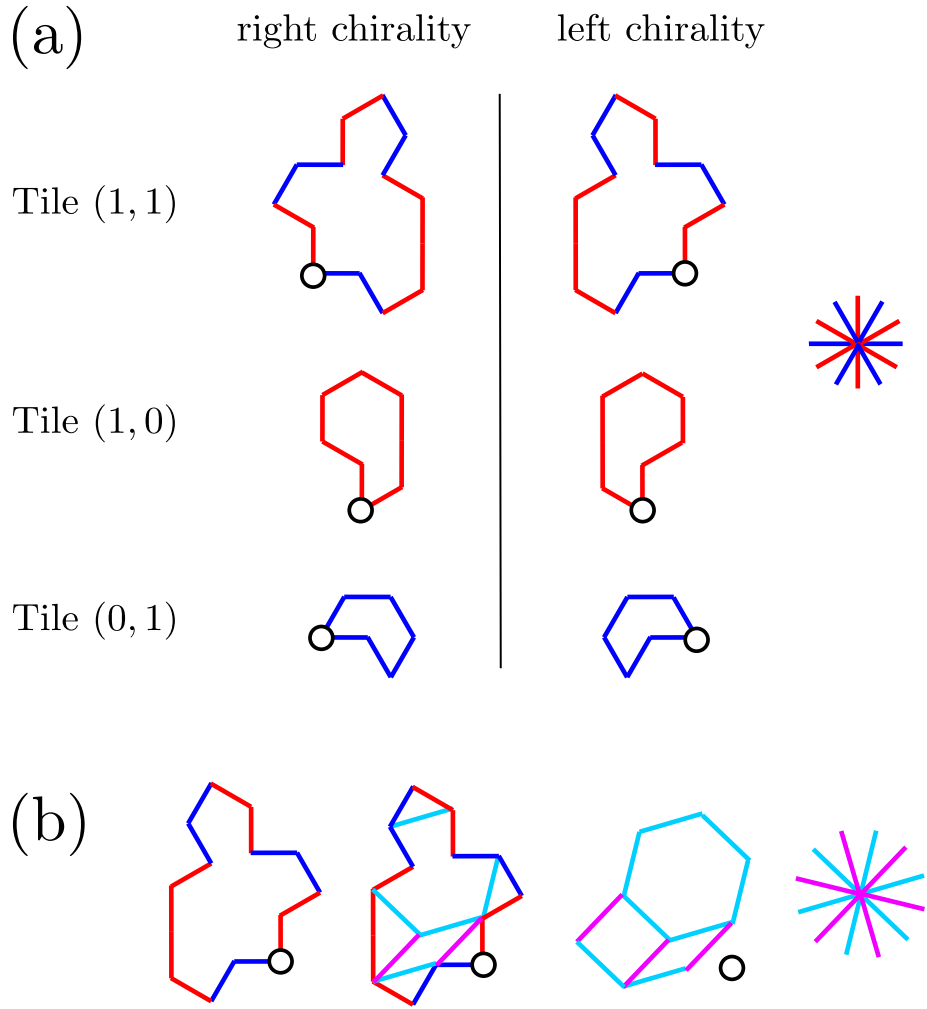}
\caption{Chiral monotile. (a) Left and right versions. Contour of Tile(1,1) is a sequence of 14 equal length edges and rotation angles of $\pi/12$. By convention, odd edges are depicted in blue and even ones in red. Tile(1,0) contains only the even edges and has a comet shape while Tile(0,1) is obtained from the blue edges and has a chevron shape. (b) Decoration of the left chiral monotile with light blue and magenta edges. Decoration edges have a length $\sqrt{2}$ larger than edges of Tile(1,1) and are rotated by $\pi/12$ with respect to them. Decoration is giving rise to one hexagon, one square and one thin rhombus. For comparison, a reference point (white circle) is depicted on all shapes.}
\label{fig:left_right_monotile}
\end{figure*}

The shape of Tile(1,1) is recalled in Figure \ref{fig:left_right_monotile} with its two chiral versions. Its contour is a specific sequence of 14 equal length edges, oriented along unit vectors of the dodecagonal symmetry. Note that one of them is repeated twice, giving rise to a long edge in the contour. Six of them (blue color) are in odd orientation and 8 of them (red color) are in even orientation. When the length of the odd edges is reduced to zero, a Tile(1,0) is obtained which has a comet shape. Conversely, Tile(0,1) is obtained when the length of the even edges is decreases to zero and has a chevron shape. By convention, we assign the left and right handed chiral forms of the monotile by looking at the two comet chiral forms. Then in all the figures of this work, we chose to always consider the left chiral form of the monotile. Equivalent derivation for the right handed monotile can be easily obtained by applying a global mirror operation. As a decoration of the chiral monotile with one hexagon, one square and one thin rhombus has been recently proposed in the literature \cite{akiyama_2025,smith_james_2024,cheritat_2024} it is also depicted in Figure \ref{fig:left_right_monotile} for comparison purposes. 
\par
\begin{figure*}
\includegraphics[width=\textwidth]{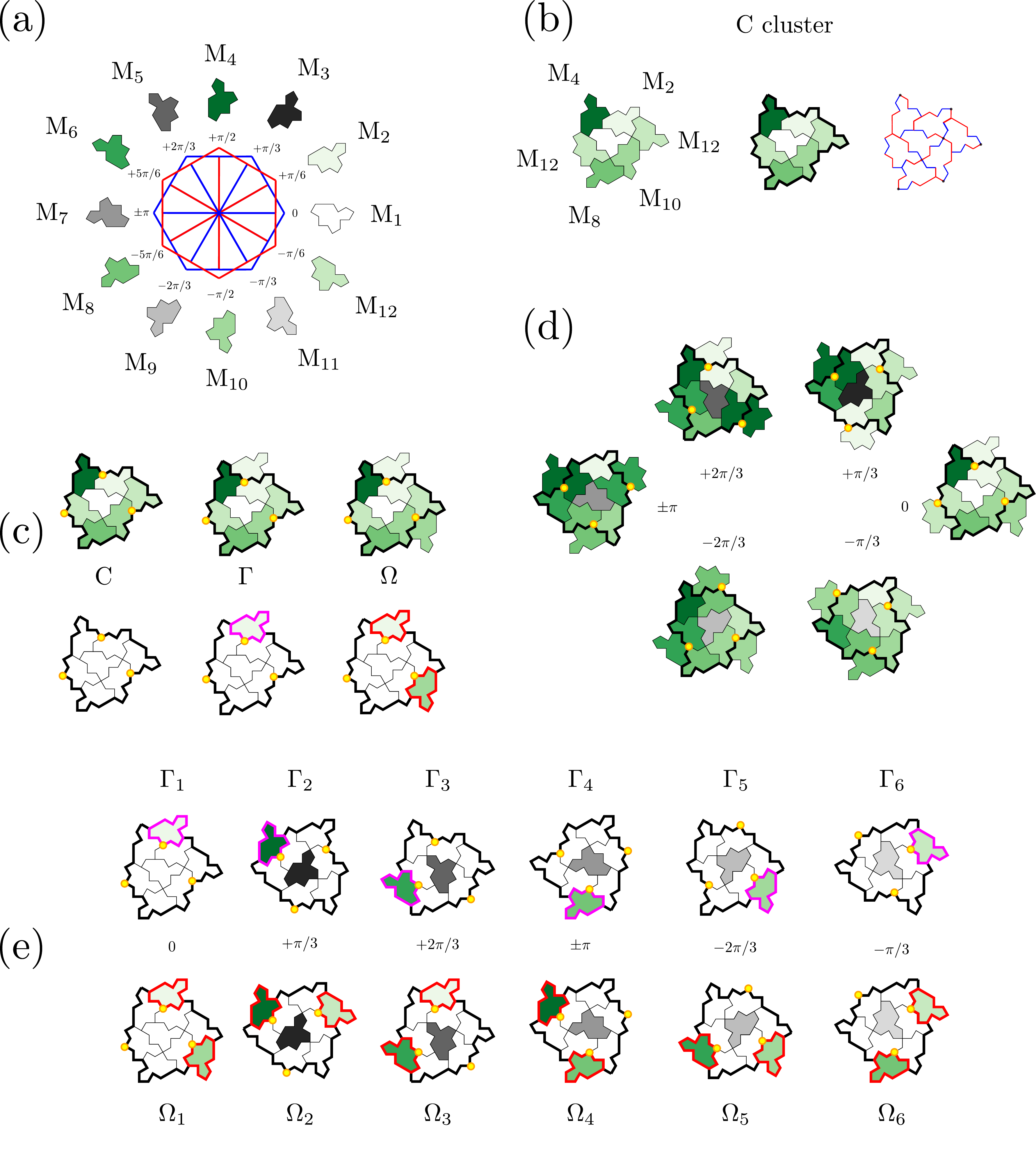}
\caption{Monotile, clusters and metatiles (a) Color code for the 12 orientations of the Tile(1,1) monotile. A left handed chirality is chosen. Shades of gray correspond to odd monotiles (M$_1$, M$_3$, ...) and shades of green to even ones (M$_2$, M$_4$, ...). (b) Construction of a C cluster of seven monotiles. A monotile (here M$_1$) is surrounded by 6 monotiles (shades of green) and the border of the C cluster is depicted by a thick black line. The same C cluster is also shown using the edge color representation (blue and red) (c) The 3 types of cluster: C, $\Gamma$ and $\Omega$. A $\Gamma$ cluster (8 monotiles) is a C cluster plus one additional even duplicated monotile (here M$_2$) while an $\Omega$ cluster (9 monotiles) contains one more additional even duplicated monotile (here M$_{10}$). (d) The six orientations of the odd monotiles in the tiling along with their nine neighboring even monotiles. Glue points (yellow color) are placed at the junction point of three even monotiles where two of them are duplicated. (e) The set of twelve metatiles corresponding to the $\Gamma$ and $\Omega$ clusters in their six orientations which allow to tile the plane.}
  \label{fig:tiles_clusters}
\end{figure*}

In a quasiperiodic tiling, this left handed monotile is taking twelve different orientations. In order to distinguish them, we introduce a color code as shown in Figure \ref{fig:tiles_clusters}(a). Its specificity is to distinguish two subsets of odd and even orientations, with two shades of color, six shades of gray (including white and black) for odd orientations, and six different shades of green for the even orientations. Shades of gray correspond to rotations by $n\pi/3$ (odd monotiles M$_1$, M$_3$, ..., M$_{11}$) and shades of green to rotations by $\pi/6+n\pi/3$ (even monotiles M$_2$, M$_4$, ..., M$_{12}$). With this convention, the repartition in the tiling of the 12 orientations of the monotile is evidenced in a convenient way. The odd orientations (shades of gray) appear for isolated monotiles surrounded by a much larger number of monotiles in even orientations (shades of green), in a ratio between 7 and 8. As computed in section \ref{sec:orientation}, the limit of this ratio is the irrational value $4+\sqrt{15} \approx 7.873$ for the infinite quasiperiodic tiling.
\par
\begin{figure*}
\includegraphics[width=\textwidth]{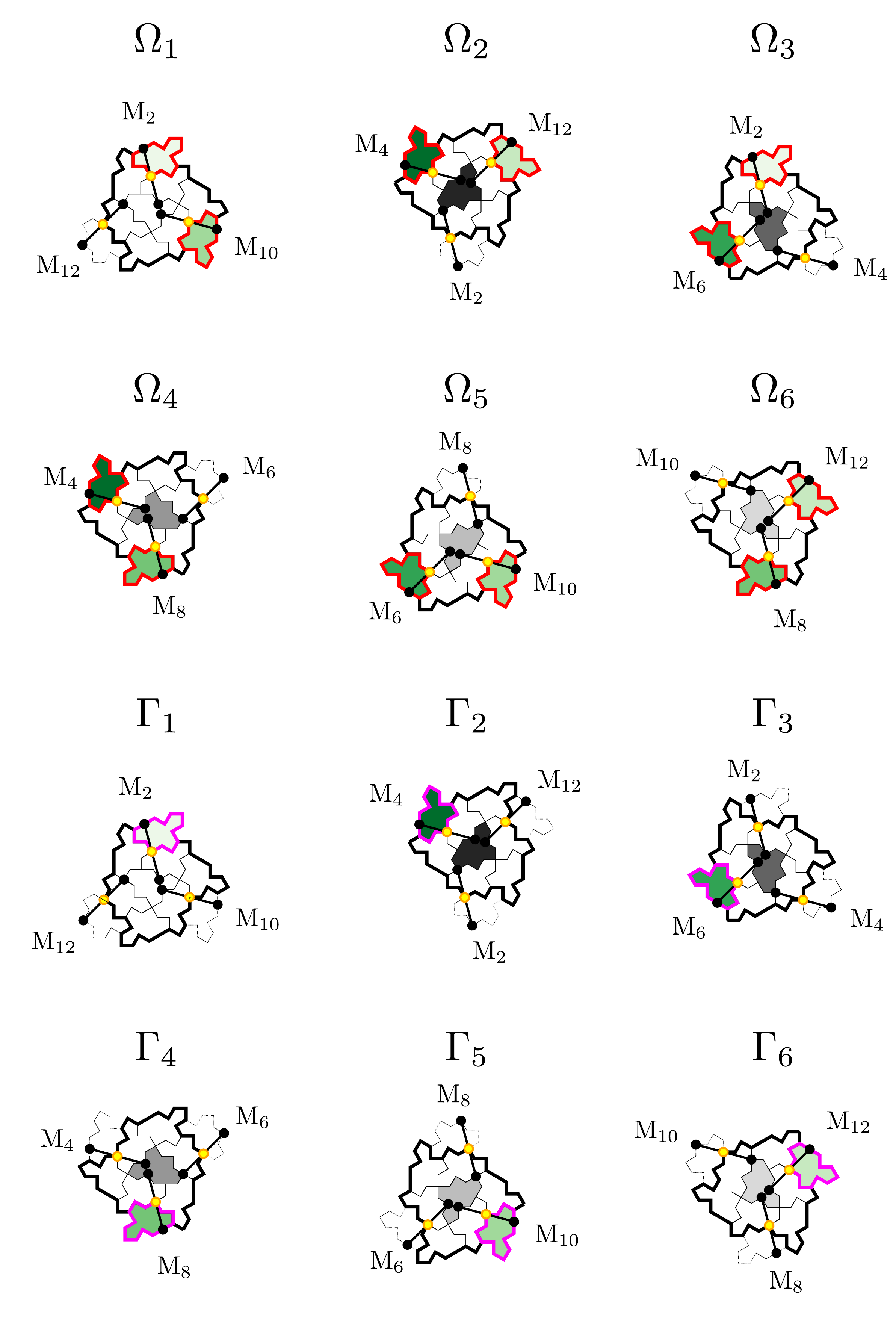}
\caption{The 12 metatiles with connection lines towards adjacent metatiles. Monotiles with a dashed line border are located in an adjacent metatile and give its position. Connection lines are connecting duplicated monotiles.}
  \label{fig:metatiles}
\end{figure*}

Each odd monotile is located at the center of a C cluster of 7 monotiles. As shown in Figure \ref{fig:tiles_clusters}(b), the odd monotile is surrounded by six even monotiles forming a corona around it. Note here that this C cluster pattern of a monotile surrounded by a corona might be related to the Heesch's problem in the field of tilings description. Such C cluster of monotiles is repeating itself in the quasiperiodic tiling by taking six different possible orientations. Doing so, the monotile at the center takes the six possible odd orientations. Moreover, another characteristic feature is that there are always three even monotiles of the corona that are duplicated by translation in the quasiperiodic tiling, as depicted in Figure \ref{fig:tiles_clusters}, and to which we can assign connection lines (Figures \ref{fig:translation} and \ref{fig:metatiles}) that will be used afterwards to determine the relative position of two adjacent metatiles (see section \ref{sec:translation}).

\section{Decomposition using two clusters}
\label{sec:decomposition}

\begin{figure*}
\includegraphics[width=\textwidth]{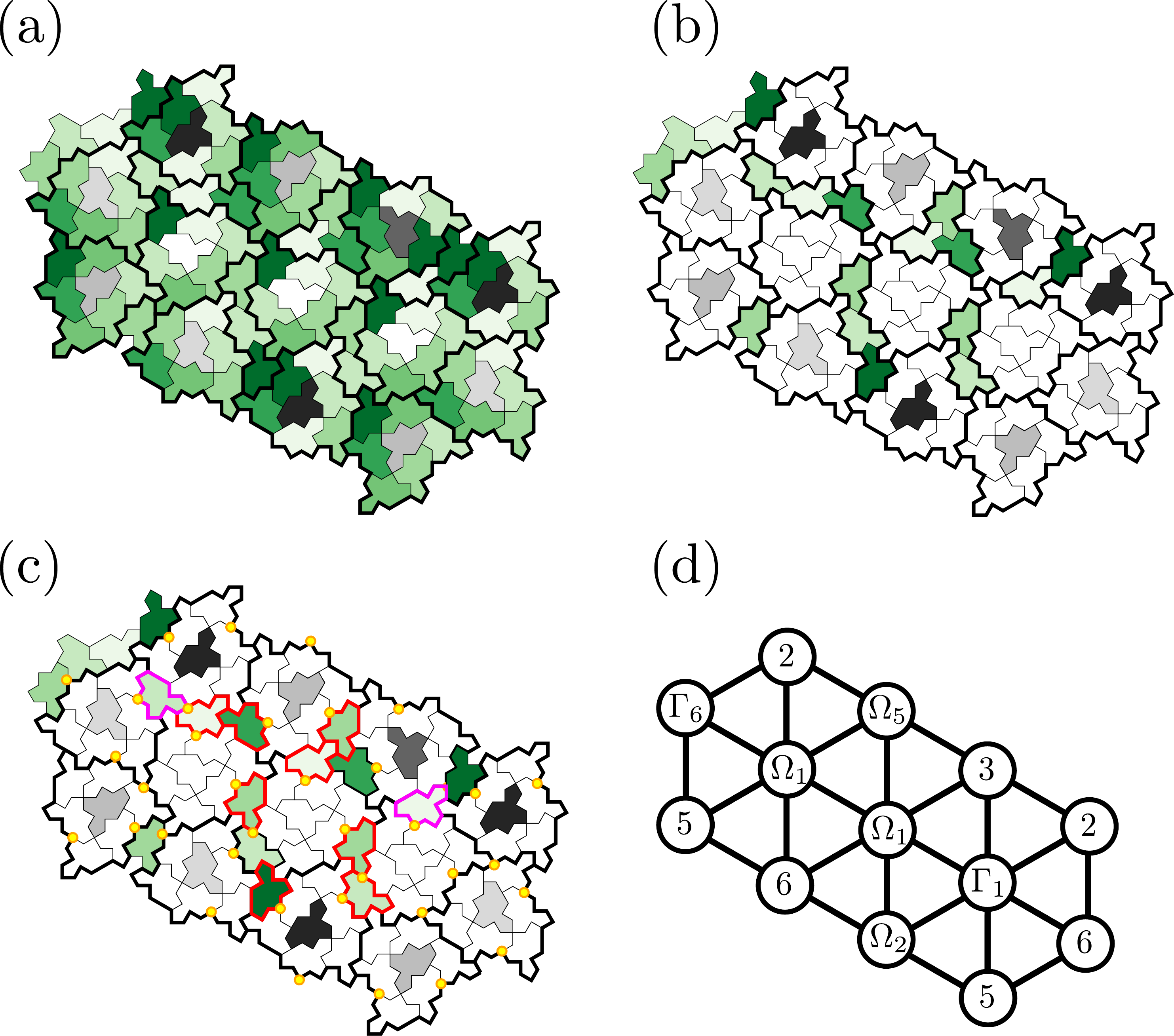}
\caption{Analysis of a typical finite region of the quasiperiodic tiling. (a) C clusters are identified by their contours (thick black lines). (b) Same region where only the orphan even monotiles outside the C clusters are evidenced (shades of green). (c)  Attribution of the orphan monotiles to either an $\Omega$ (orange contour) or a $\Gamma$ (pink contour) metatile. (d) Mapping onto a triangular lattice of the metatiles. Attribution to a metatile needs the knowledge of enough adjacent clusters. Otherwise, at the border of the finite region, the number corresponds to the orientation of the metatile, as it might be either a $\Gamma$ or an $\Omega$ one.}
\label{fig:example}
\end{figure*}

Indeed, in order to describe the tiling, we introduce a decomposition into only two clusters noted $\Gamma$ and $\Omega$, which contain respectively eight and nine monotiles. To illustrate this decomposition, let's consider a typical finite region of a quasiperiodic tiling depicted in Figure \ref{fig:example}(a). This region was obtained using the applet available at \url{https://cs.uwaterloo.ca/~csk/spectre/app.html} developed by the authors of \cite{smith2023chiral}. One can identify the C clusters by their contour and orientation. Doing so, we evidence the fact that most of the monotiles are inside C clusters, expect some extra orphan monotiles in even orientation as shown in Figure \ref{fig:example}(c).
\par
The main idea is to attribute each orphan even monotile to a nearby C cluster in a unique way. This is done by introducing the $\Gamma$ and $\Omega$ clusters as shown in Figure \ref{fig:tiles_clusters}(c). Note these two types of clusters are already introduced in \cite{smith2023chiral}. These two clusters are indeed C clusters to which one ($\Gamma$) or two ($\Omega$) orphan even monotiles are attached. An intricate part is assigning the orphan even monotiles to the C clusters in a unique way. When two orphan monotiles are 'attached' to a C cluster, a $\Omega$ cluster is identified. When only one orphan monotile is attached, the cluster type is a $\Gamma$ one. The attribution can be done manually in a unique way for the finite region chosen as an example as shown in Figure \ref{fig:example}(e). 
\par
As a result, the decomposition of the entire tiling contains only two clusters $\Gamma$ and $\Omega$. Because they are taking six possible orientations, we end up with a puzzle of twelve metatiles which can be respectively labeled as ($\Gamma_1$, $\Gamma_2$, $\Gamma_3$, $\Gamma_4$, $\Gamma_5$, $\Gamma_6$) and ($\Omega_1$, $\Omega_2$, $\Omega_3$, $\Omega_4$, $\Omega_5$, $\Omega_6$) as shown in Figures \ref{fig:tiles_clusters}(e) and \ref{fig:metatiles}. The index from 1 to 6 refers to the orientation of the odd monotile at the center of each metatile. Compared to the original publication \cite{smith2023chiral}, the clusters are essentially the same, but their six possible orientations are explicitly encoded by the index. The uniqueness of this decomposition into twelve metatiles is demonstrated by the inflation construction.
\par

\par
\section{A puzzle of twelve metatiles}
\label{sec:translation}
Using twelve metatiles allow to properly separate translation and rotation operations. Indeed, as all the possible orientations are already incorporated in the metatiles, only translations are needed between metatiles to construct the quasiperiodic tiling.
\par
\begin{figure*}
\includegraphics[width=\textwidth]{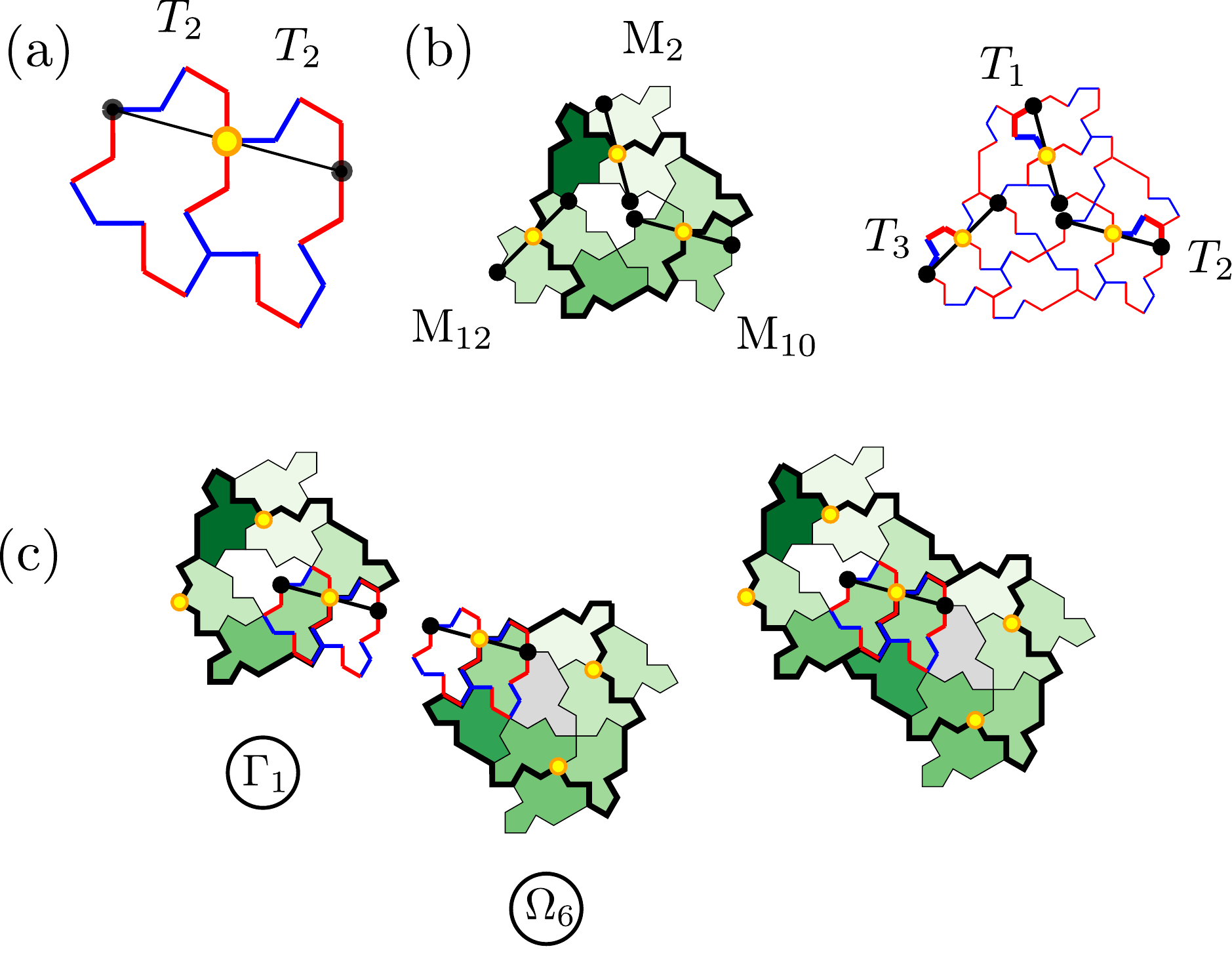}
\caption{Translation between two adjacent metatiles. (a) Translation (here $T_2$) between two adjacent monotiles in the same orientation (here M$_{10}$). A glue point (yellow color) is placed as a marker. Translation to superimpose the two monotiles is from a yellow to a black point and conversely. (b) In the quasiperiodic tiling, around each odd monotile (here M$_1$) three even monotiles (here M$_{10}$, M$_{12}$ and M$_{2}$) are always duplicated twice. The three corresponding glue points (yellow color) are located on the border of the C cluster (thick back line). Translation vectors between two adjacent monotiles are $T_1$ (M$_2$/M$_8$), $T_2$ (M$_4$/M$_{10}$) and $T_3$ (M$_6$/M$_{12}$) with $T_3=T_1+T_2$. (c) Example with translation $T_2$ between $\Gamma_1$ and $\Omega_6$ metatiles with two M$_{10}$ duplicated monotiles.}
  \label{fig:translation}
\end{figure*}

However, determining all the translations between two adjacent metatiles is not an easy task. However, we propose to take advantage of the fact that, as previously mentioned in section \ref{sec:tiles_clusters}, some even monotiles are duplicated with each of them belonging to two different but adjacent metatiles. As shown in Figure \ref{fig:translation}, we can introduce connection lines which length is twice the translation between two duplicated monotiles. An associated glue point (yellow color) is placed as a reference point which is located on the border of the C cluster (black thick line). 
\par
In the whole tiling, the three directions along which adjacent even monotiles are duplicated deduce from each other by a $2\pi/3$ rotation, giving rise to three translation vectors $T_1$, $T_2$ and $T_3=T_1+T_2$ between duplicated even monotiles. Connection lines are taking the same three directions. It is convenient to code edge directions of the monotile using the 12 unit length vectors deduced from each other by $2\pi/12$ rotation using the 4 vectors basis $e_1$, $e_2$, $e_3$ and $e_4$. The 12 unit vectors read: $e_1$, $e_2$, $e_3$, $e_4$, $e_3-e_1$, $e_4-e_2$, $-e_1$, $-e_2$, $-e_3$, $-e_4$, $e_1-e_3$, $e_2-e_4$. Advantage is that this four vector basis is already containing the lifted version in a four dimensional space \cite{imperor-clerc_square-triangle_2021,Imperor_PRB_2024}. Moreover, by changing the ratio between the lengths of odd and even unit vectors, tilings made of Tile (a,b) instead of Tile(1,1) are obtained directly. Expressed in this 4 vector basis, the translation vectors $T_1$, $T_2$ and $T_3$ are: 
\begin{enumerate}
    \item[-] M$_2$/M$_8$ monotiles: $T_1=-2e_1+e_2+e_3+e_4$
    \item[-] M$_4$/M$_{10}$ monotiles: $T_2=e_1+e_2+e_3-2e_4$
    \item[-] M$_6$/M$_{12}$ monotiles: $T_3=T_1+T_2=-e_1+2e_2+2e_3-e_4$
\label{eq:translations}
\end{enumerate}
\par
\begin{figure*}
\includegraphics[width=1.2\textwidth]{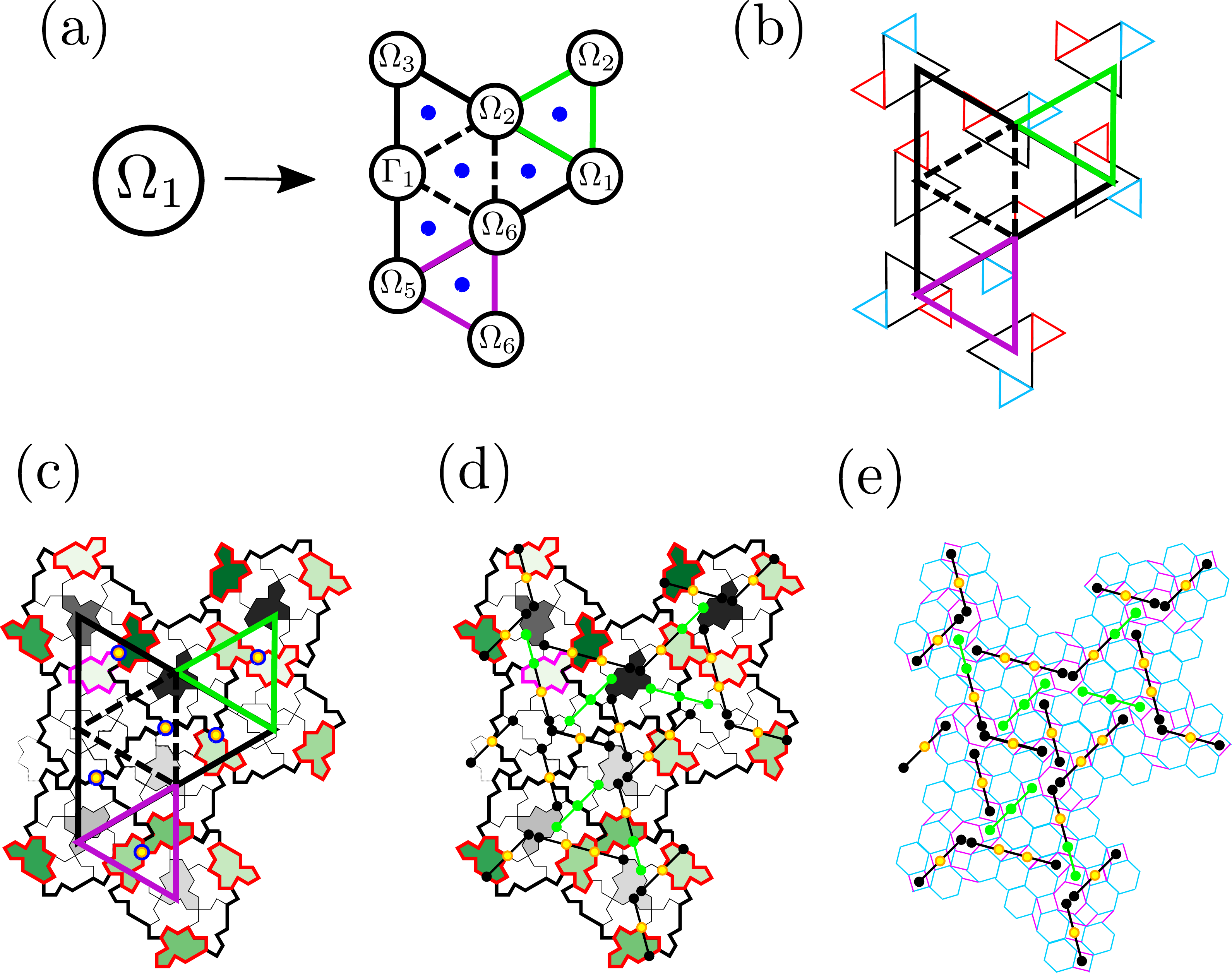}
\caption{First inflation step. Substitution rule for the $\Omega_1$ metatile mapped on a triangular lattice (a). The corresponding triangular patterns (b). Assembly of metatiles superimposed with the triangular pattern (c), the connection lines and glue points (d) and the decoration representation (e), where glue points are at the center of some squares. Green junction lines are placed at the center of the remaining squares of the decoration and correspond to pairs of adjacent even monotiles turned by $\pi$. For the first inflation step, junction points between three metatiles (dark blue) superimpose on glue points (yellow). }
  \label{fig:inflation_step1}
\end{figure*}

The twelve metatiles as shown in Figure \ref{fig:metatiles} along with the connection lines between even monotiles, where a duplicated monotile belonging to an adjacent metatile is depicted with a dashed line border. In this way, the translation is known for some pairs of adjacent metatiles, allowing to place them with respect to each other like illustrated in Figure \ref{fig:translation}. A characteristic feature of the quasiperiodic tiling is that duplicated even monotiles always appear in groups of two or three, as can be visualized by the connection lines, that carry either one or two yellow glue points as shown for example in Figure \ref{fig:inflation_step1}.
\par
In addition, green junction lines can be depicted to complete the pattern of junction lines. Remarkably, a green junction line corresponds to a pair of adjacent even monotiles turned by a rotation angle $\pi$ and translated, as can be seen for example in Figure \ref{fig:inflation_step1}(d). The length and the three possible orientations of a green junction line are again given by $T_1$, $T_2$ and $T_3$. As discussed furthermore in section \ref{sec:dual}, all junction lines are forming triangular patterns of different types, in accordance with the mapping of the metatiles assembly onto a triangular lattice.

\section{Mapping on a triangular lattice and inflation}
\label{sec:inflation}
To describe the inflation construction, we rely on the fact that an assembly of twelve metatiles can be mapped onto a triangular lattice, as shown on a finite size region in Figure \ref{fig:example}(d). This mapping is containing all the information about how the metatiles are arranged together.
\par
\begin{figure*}
\includegraphics[width=1.2\textwidth]{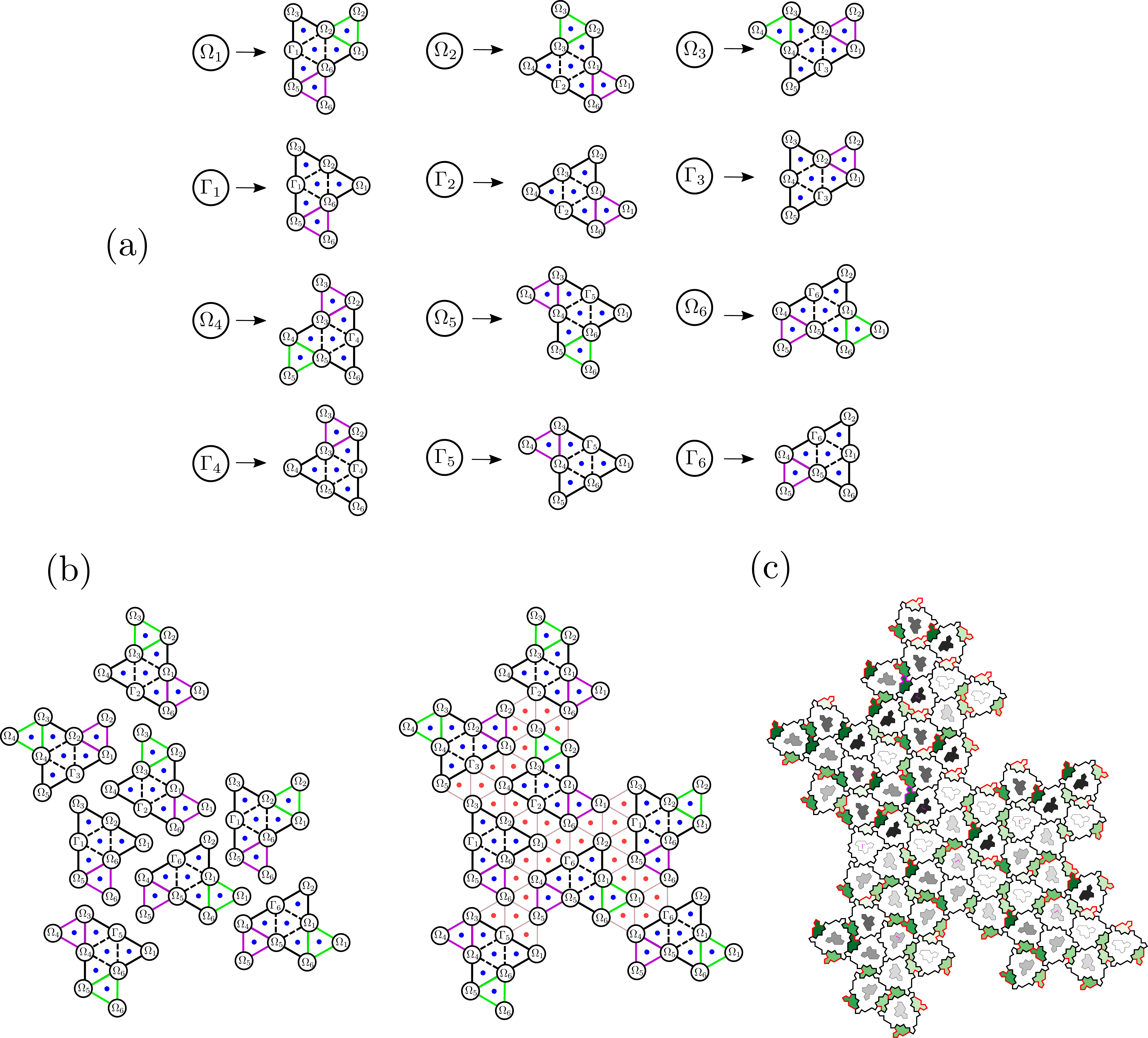}
\caption{Second inflation step. (a) Substitution rule for the 12 metatiles mapped on a triangular lattice. Junction points between three metatiles are depicted by dark blue dots. (b) Second inflation step mapped on a triangular lattice using $\Omega_1$ metatile as a seed. New junction points between metatiles are depicted by red dots. (c) The corresponding assembly of metatiles.}
  \label{fig:inflation_step2}
\end{figure*}

As a starting seed for the inflation construction, either a $\Gamma$ or an $\Omega$ cluster can be chosen. For convenience, the $\Omega$ cluster in its first orientation, corresponding to the $\Omega_1$ metatile, is considered as the seed in the following Figures \ref{fig:omega_gamma}, \ref{fig:inflation_step1} and \ref{fig:inflation_step2}. Starting from the seed, substitution rules are needed for constructing the quasiperiodic tiling. The substitution rule for the two clusters $\Omega$ and $\Gamma$ is given in Figure \ref{fig:omega_gamma}. It is adapted from \cite{smith2023chiral} by taking explicitly into account the orientation of the monotiles. Moreover, the chirality of the monotile is always fixed, in contrast with the alternance in chiralty at each inflation step used in \cite{smith2023chiral}. First, we identify in an $\Omega$ cluster a triangular pattern of six triangles with eight monotiles placed on its vertices. This triangular pattern is chiral as identified by the red/blue colors used for the two triangles in the pattern corresponding to duplicated monotiles. Note that the triangular pattern for a $\Gamma$ cluster is derived in the same way by omitting the blue triangle. 
\par
The substitution rule is detailed in the caption of Figure \ref{fig:omega_gamma} and is based on a relationship between a triangular pattern and the orientations of the monotiles inside it. In the substitution, an $\Omega$ metatile is replaced by 8 metatiles (one of type $\Gamma$ and 7 of type $\Omega$) when a $\Gamma$ metatile is replaced by 7 metatiles (one of type $\Gamma$ and 6 of type $\Omega$). Note that the substitution rule is essentially the same for both metatiles, except that one $\Omega$ metatile is omitted. The first inflation step is shown in Figure \ref{fig:inflation_step1}. 
\par
\begin{figure*}
\includegraphics[width=\textwidth]{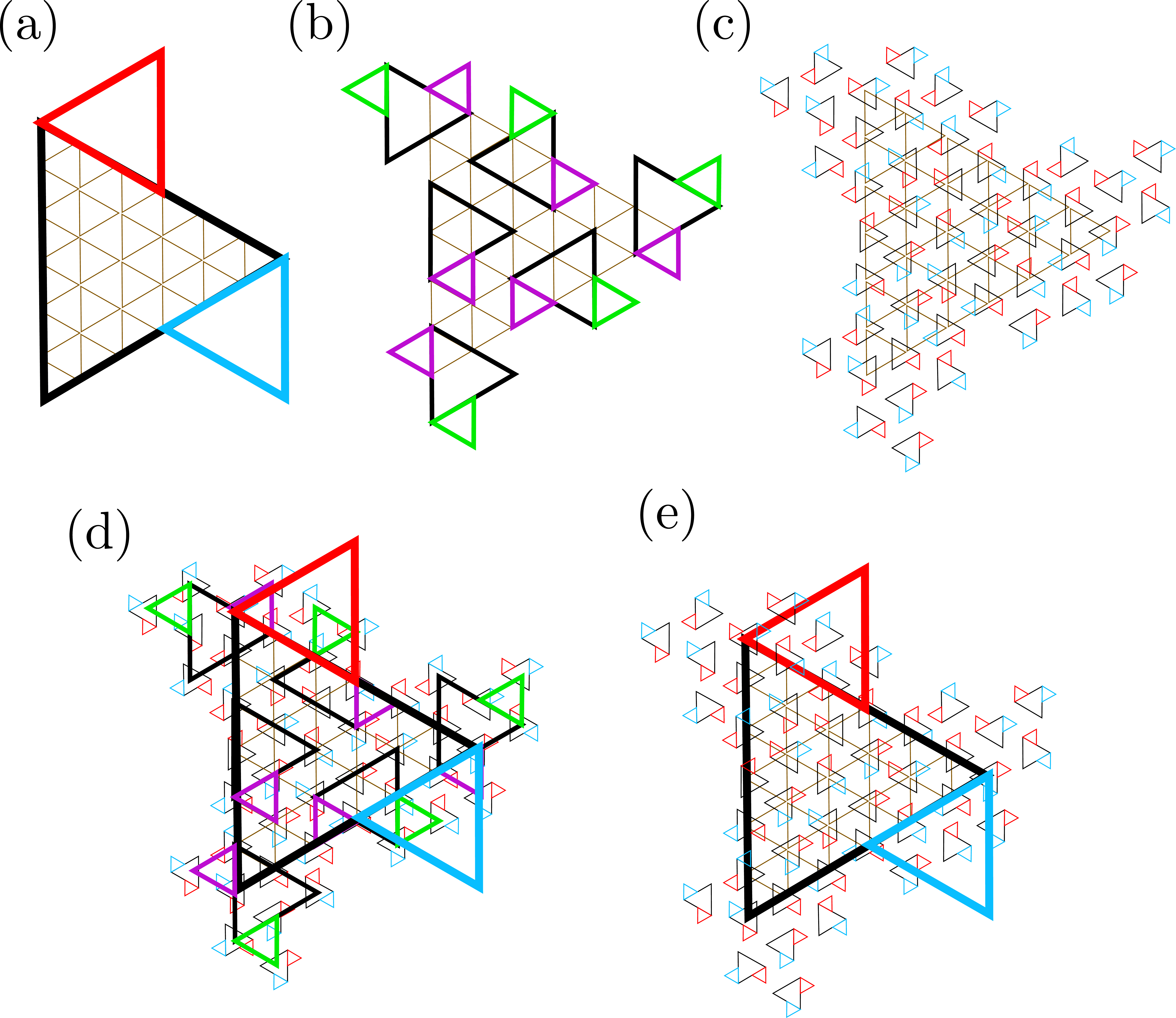}
\caption{Comparison of the two first inflation steps showing the alternance in chirality of the triangular patterns. (a) Triangular pattern at step 2 with blue/red triangles. (b) Triangular patterns at step 1 with green/magenta triangles. (c) Triangular patterns associated to the metatiles (red/blue triangles). (d) Superposition of all triangular patterns. (e) Same superposition without step 1: chirality of the triangular patterns is the same for the metatiles and at step 2.}
  \label{fig:inflation}
\end{figure*}

An important feature of the homochiral substitution rule is that it is inducing an alternance in chirality of the triangular patterns depending on the number of inflation steps, as illustrated in Figure \ref{fig:inflation}. Indeed, after two inflation steps, the initial chirality of the triangular pattern (depicted in red and blue colors) is recovered, when the opposite chirality at step 1 is obtained (depicted in green and magenta colors). This alternance in chirality can be understood as the signature of chirality when using an homochiral inflation construction. The resulting tiling is the same as in \cite{smith2023chiral}, as can be seen for example in Figure \ref{fig:inflation_step2}(c).   
\par

\section{Orientation distribution of the monotile}
\label{sec:orientation}
From the substitution rule given in Figure \ref{fig:omega_gamma}, the total numbers of $\Gamma$ and $\Omega$ metatiles at each inflation step can be computed using a 2x2 inflation matrix:
\begin{equation}
    \begin{pmatrix}
        N_{\Gamma}\\N_{\Omega}
    \end{pmatrix}
    \mapsto
    \begin{pmatrix}
        1 & 1\\
        6 & 7
    \end{pmatrix}
    \begin{pmatrix}
        N_{\Gamma}\\N_{\Omega}
    \end{pmatrix}
    \label{eq:2by2matrix}
\end{equation}
For example, starting from the $\Gamma_1$ metatile as a seed, one obtain the following number of metatiles for the first inflation steps:
\begin{enumerate}
    \item    $N_{\Gamma}$=1   $N_{\Omega}$=6  $N_{\Omega}$/$N_{\Gamma}$=6
    \item    $N_{\Gamma}$=7   $N_{\Omega}$=48  $N_{\Omega}$=6  $N_{\Omega}$/$N_{\Gamma}$=6.8571429
    \item    $N_{\Gamma}$=55   $N_{\Omega}$=378  $N_{\Omega}$=6  $N_{\Omega}$/$N_{\Gamma}$=6.8727273
    \item    $N_{\Gamma}$=433  $N_{\Omega}$=2976   $N_{\Omega}$/$N_{\Gamma}$=6.8729792
    \item    $N_{\Gamma}$=3409   $N_{\Omega}$=23430  $N_{\Omega}$/$N_{\Gamma}$=6.8729833
\end{enumerate}

The maximum eigenvalue, also called Perron root, of the 2x2 inflation matrix is $4+\sqrt{15} \approx 7.873$, and the associated eigenvector gives the proportion of the two types of metatiles in the infinite quasiperiodic tiling:
\begin{equation}
  N^{\infty}_{\Omega}=(3+\sqrt{15}) N^{\infty}_{\Gamma} \approx 6.8729833~N^{\infty}_{\Gamma}
\label{eq:proportion_omega_gamma}
\end{equation}
Note that the convergence of the inflation is rather fast, with a ratio $N_{\Omega}$/$N_{\Gamma}$ already very close to the limit value after only 5 inflation steps.
\par
Including orientation, the same substitution rule for all metatiles can be written in the form of a 12x12 matrix $M_{\Gamma,\Omega}$:
\begin{equation}
    \begin{pmatrix}
        N_{\Gamma_1}\\ N_{\Gamma_2}\\
        N_{\Gamma_3}\\ N_{\Gamma_4}\\
        N_{\Gamma_5}\\ N_{\Gamma_6}\\
        N_{\Omega_1}\\
        N_{\Omega_2}\\
        N_{\Omega_3}\\
        N_{\Omega_4}\\
        N_{\Omega_5}\\
        N_{\Omega_6}\\
    \end{pmatrix}
    \mapsto
    \begin{pmatrix}
    1, 0, 0, 0, 0, 0, 1, 0, 0, 0, 0, 0\\
    0, 1, 0, 0, 0, 0, 0, 1, 0, 0, 0, 0\\
    0, 0, 1, 0, 0, 0, 0, 0, 1, 0, 0, 0\\
    0, 0, 0, 1, 0, 0, 0, 0, 0, 1, 0, 0\\
    0, 0, 0, 0, 1, 0, 0, 0, 0, 0, 1, 0\\
    0, 0, 0, 0, 0, 1, 0, 0, 0, 0, 0, 1\\
    1, 2, 1, 0, 1, 1, 1, 2, 1, 0, 1, 2\\
    1, 1, 2, 1, 0, 1, 2, 1, 2, 1, 0, 1\\
    1, 1, 1, 2, 1, 0, 1, 2, 1, 2, 1, 0\\
    0, 1, 1, 1, 2, 1, 0, 1, 2, 1, 2, 1\\
    1, 0, 1, 1, 1, 2, 1, 0, 1, 2, 1, 2\\
    2, 1, 0, 1, 1, 1, 2, 1, 0, 1, 2, 1
    \end{pmatrix}
    \begin{pmatrix}
        N_{\Gamma_1}\\ N_{\Gamma_2}\\
        N_{\Gamma_3}\\ N_{\Gamma_4}\\
        N_{\Gamma_5}\\ N_{\Gamma_6}\\
        N_{\Omega_1}\\
        N_{\Omega_2}\\
        N_{\Omega_3}\\
        N_{\Omega_4}\\
        N_{\Omega_5}\\
        N_{\Omega_6}\\
    \end{pmatrix}
    \label{eq:12by12matrix_omega_gamma}
\end{equation}
This 12x12 $M_{\Gamma,\Omega}$ matrix is expressing the same substitution rule as illustrated in Figure \ref{fig:omega_gamma} but for all six orientations of the $\Gamma$ and $\Omega$ clusters giving rise to the twelve metatiles. Amounts of all metatiles at successive inflation steps are computed from successive powers of the matrix $M_{\Gamma,\Omega}$.
\par
Finally, from the matrix $M_{\Gamma,\Omega}$, one can compute the orientation distribution on the monotile in terms of its 12 orientations ($M_1$, $M_2$ ... $M_{12}$) for successive inflation steps. First, the composition of all metatiles in terms of the 12 orientations of the monotile is written as a 12x12 matrix $M$:
\begin{equation}
    \begin{pmatrix}
        N_{M_1}\\ N_{M_2}\\
        N_{M_3}\\ N_{M_4}\\
        N_{M_5}\\ N_{M_6}\\
        N_{M_7}\\
        N_{M_8}\\
        N_{M_9}\\
        N_{M_{10}}\\
        N_{M_{11}}\\
        N_{M_{12}}\\
    \end{pmatrix}
    =
    \begin{pmatrix}
    1, 0, 0, 0, 0, 0, 1, 0, 0, 0, 0, 0\\
    2, 2, 1, 1, 0, 1, 2, 2, 2, 1, 0, 1\\
    0, 1, 0, 0, 0, 0, 0, 1, 0, 0, 0, 0\\
    1, 2, 2, 1, 1, 0, 1, 2, 2, 2, 1, 0\\
    0, 0, 1, 0, 0, 0, 0, 0, 1, 0, 0, 0\\
    0, 1, 2, 2, 1, 1, 0, 1, 2, 2, 2, 1\\
    0, 0, 0, 1, 0, 0, 0, 0, 0, 1, 0, 0\\
    1, 0, 1, 2, 2, 1, 1, 0, 1, 2, 2, 2\\
    0, 0, 0, 0, 1, 0, 0, 0, 0, 0, 1, 0\\
    1, 1, 0, 1, 2, 2, 2, 1, 0, 1, 2, 2\\
    0, 0, 0, 0, 0, 1, 0, 0, 0, 0, 0, 1\\
    2, 1, 1, 0, 1, 2, 2, 2, 1, 0, 1, 2
    \end{pmatrix}
    \begin{pmatrix}
        N_{\Gamma_1}\\ N_{\Gamma_2}\\
        N_{\Gamma_3}\\ N_{\Gamma_4}\\
        N_{\Gamma_5}\\ N_{\Gamma_6}\\
        N_{\Omega_1}\\
        N_{\Omega_2}\\
        N_{\Omega_3}\\
        N_{\Omega_4}\\
        N_{\Omega_5}\\
        N_{\Omega_6}\\
    \end{pmatrix}
    \label{eq:m_metatiles_mono}
\end{equation}
Then the inflation matrix $M_{mono}$ in terms of the 12 orientations of the monotile equals:
\begin{equation}
  M_{mono}=M \cdot M_{\Gamma,\Omega} \cdot M^{-1}
\label{eq:matrix_mono_calculation}
\end{equation}
and one gets the coefficients of the matrix $M_{mono}$:
\begin{equation}
    \begin{pmatrix}
        N_{M_1}\\ N_{M_2}\\
        N_{M_3}\\ N_{M_4}\\
        N_{M_5}\\ N_{M_6}\\
        N_{M_7}\\
        N_{M_8}\\
        N_{M_9}\\
        N_{M_{10}}\\
        N_{M_{11}}\\
        N_{M_{12}}\\
    \end{pmatrix}
    \mapsto
    \begin{pmatrix}
    1, 0, 2, 0, 0, 0, -2, 1, -1, 0, 0, 0\\
    0, 1, 5, 0, 3, 1, 0, 2, -4, 2, -5, 2\\
    0, 0, 1, 0, 2, 0, 0, 0, -2, 1, -1, 0\\
    -5, 2, 0, 1, 5, 0, 3, 1, 0, 2, -4, 2\\
    -1, 0, 0, 0, 1, 0, 2, 0, 0, 0, -2, 1\\
    -4, 2, -5, 2, 0, 1, 5, 0, 3, 1, 0, 2\\
    -2, 1, -1, 0, 0, 0, 1, 0, 2, 0, 0, 0\\
    0, 2, -4, 2, -5, 2, 0, 1, 5, 0, 3, 1\\
    0, 0, -2, 1, -1, 0, 0, 0, 1, 0, 2, 0\\
    3, 1, 0, 2, -4, 2, -5, 2, 0, 1, 5, 0\\
    2, 0, 0, 0, -2, 1, -1, 0, 0, 0, 1, 0\\
    5, 0, 3, 1, 0, 2, -4, 2, -5, 2, 0, 1
    \end{pmatrix}
    \begin{pmatrix}
        N_{M_1}\\ N_{M_2}\\
        N_{M_3}\\ N_{M_4}\\
        N_{M_5}\\ N_{M_6}\\
        N_{M_7}\\
        N_{M_8}\\
        N_{M_9}\\
        N_{M_{10}}\\
        N_{M_{11}}\\
        N_{M_{12}}\\
    \end{pmatrix}
    \label{eq:12by12matrix_monotile}
\end{equation}
Starting from any metatile as a seed, one can calculate the values of ($N_{M_1}, ... ,N_{M_{12}}$) for all inflation steps. For example, by taking metatile $\Gamma_1$ as a seed, the initial repartition is:
\begin{equation*}
(N_{M_1}, ... ,N_{M_{12}})=(1, 2, 0, 1, 0, 0, 0, 1, 0, 1, 0, 2)
\end{equation*}
and the numbers of monotiles for the first six inflation steps are:
\begin{enumerate}
     \item    {(2, 10, 1, 7, 1, 7, 0, 9, 1, 10, 2, 12)}
    \item    {(12, 73, 9, 67, 7, 65, 7, 71, 9, 78, 11, 79)}
    \item    {(78, 569, 72, 550, 66, 550, 65, 567, 73, 585, 79, 588)}
    \item    {(586, 4470, 564, 4420, 548, 4421, 552, 4477, 571, 4528, 588, 4523)}
    \item    {(4516, 35187, 4458, 35053, 4413, 35084, 4428, 35245, 4491, 35380, 
4533, 35354)}
\item  {(35330, 277129, 35149, 276777, 35041, 276911, 35107, 277404, 35281, 
277752, 35395, 277612)}
\end{enumerate}
The maximum eigenvalue of $M_{mono}$ is $4+\sqrt{15}$ and, from its corresponding eigenvector, $(4 - \sqrt{15}, 1, 4 - \sqrt{15}, 1, 4 - \sqrt{15}, 1, 4 - \sqrt{15}, 1, 4 \
- \sqrt{15}, 1, 4 - \sqrt{15}, 1)$, we obtain the fact that in the infinite quasiperiodic tiling, the six odd monotiles (shades of gray) are in the same amount, with $N^{\infty}_{M_1}=N^{\infty}_{M_3}=N^{\infty}_{M_5}=N^{\infty}_{M_7}=N^{\infty}_{M_9}=N^{\infty}_{M_{11}}$, as well as the even monotiles (shades of green) with $N^{\infty}_{M_2}=N^{\infty}_{M_4}=N^{\infty}_{M_6}=N^{\infty}_{M_8}=N^{\infty}_{M_{10}}=N^{\infty}_{M_{12}}$. Moreover, the proportion between odd and even monotiles $N^{\infty}_{even}/N^{\infty}_{odd}$ in the infinite tiling is computed and is equal to $4+\sqrt{15} \approx 7.873$:
\begin{equation}
  N^{\infty}_{even}=(4+\sqrt{15}) N^{\infty}_{odd}
\label{eq:proportion_even_odd}
\end{equation}
a value in agreement with literature \cite{smith2023chiral}.

\section{Link with dual representation}
\label{sec:dual}
\begin{figure*}
\includegraphics[width=0.9\textwidth]{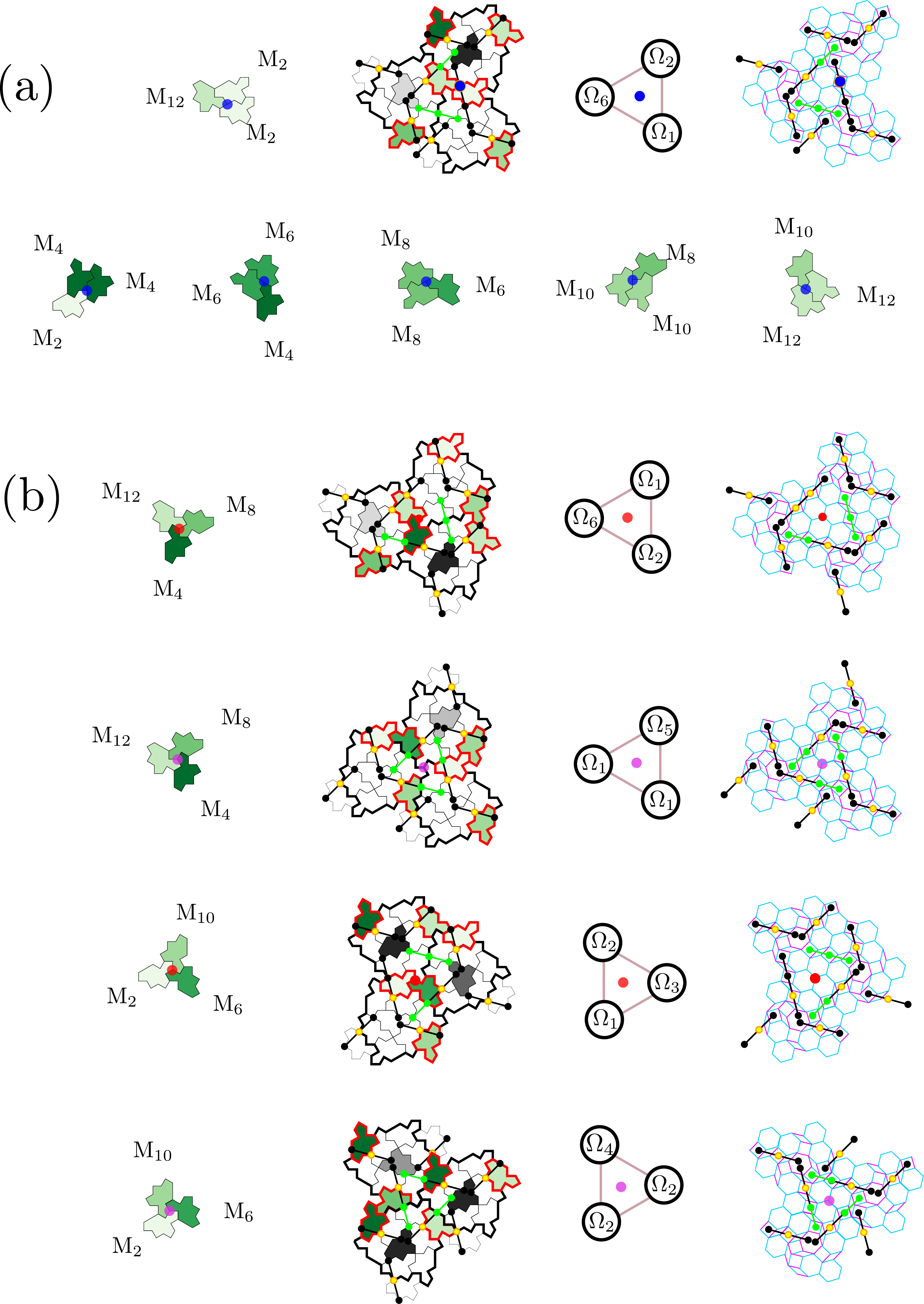}
\caption{Junction points between three adjacent metatiles. (a) First type of junction points (here in dark blue) with six possible orientations. In the example, the three monotiles at the junction point are two duplicated M$_2$ monotiles plus an M$_{12}$ monotile. Such junction points are generated during the first inflation step and superimpose on yellow glue points (see Figure \ref{fig:inflation_step1}). (b) Second type of junction points (in red and magenta) with four possible configurations. At such junction points, monotiles can be either (M$_2$, M$_6$, M$_{10}$) or (M$_4$, M$_8$, M$_{12}$). These junction points are generated during the second inflation step (see Figure \ref{fig:inflation_step2}).}
  \label{fig:junction_points}
\end{figure*}

An equivalent way to describe how the twelve monotiles are assembled together is to consider for the mapping, instead of the triangular lattice, its dual hexagonal lattice. This can be done in different ways. One possibility is to consider the set of junction points between three metatiles as illustrated in Figure \ref{fig:junction_points}. Such points are located at the intersection of three adjacent metatiles, and, based on the two first inflation steps, two types of such junction points are observed. For the first type (dark blue color), two duplicated monotiles are involved at the junction point, for example two M$_2$ monotiles with an M$_{12}$ monotile, giving rise to six different orientations for such junction points. For the second type (red and magenta colors) of junction points, a set of three different monotiles is involved, giving rise to four different combinations, involving either (M$_2$, M$_6$, M$_{10}$) or (M$_4$, M$_8$, M$_{12}$) monotiles.
\par
Another way is to use the decoration of the monotile shown in Figure \ref{fig:left_right_monotile}(b). Doing so, the connection lines (black and green colors together) superimposed on the rows of alternate squares and rhombuses, with the extremities of all junction lines located inside the hexagons (pink color) corresponding to the odd monotiles. The connection lines are forming triangles of three different sizes that are surrounding clusters of either one, three or six hexagons (light blue color) as can be seen for example in Figure \ref{fig:inflation_step1}(e). Such clusters of hexagons appear to be the building units of the dual mapping onto an hexagonal lattice \cite{smith_james_2024,cheritat_2024}. For example, in Figure \ref{fig:junction_points}(b) clusters of three or six hexagons in two possible orientations are delimitated by three junction lines of different sorts.

\section{Conclusion}
Determining the combinatorial properties of the Tile(1,1) chiral monotile can be approached in different manners. In this work, we suggest that it can be done when fixing the chirality of the monotile allowing to calculate the orientation distribution of its twelve possible orientations at different inflation steps. Note that the inflation rule we are proposing is leading to a single quasicrystalline tiling in a deterministic way. A perspective of this work would be to generate figures using a computer code for further inflation steps as so far, only the two first inflation steps have been constructed manually in all the figures. Variants of this quasiperiodic tiling should be also included in a more general description.
\par
Another important conclusion of this work is to provide different equivalent sets of objects for further investigations. One set is the twelve metatiles and their mapping onto a triangular lattice. Another set is the different types of junction lines and the different triangles they are forming directly on the tiling. This later set of triangles is most probably equivalent to the set of hexagon clusters considered in \cite{cheritat_2024} when using the decoration of the monotile by hexagons, squares and rhombuses. In this spirit, junction points between three metatiles might be an interesting alternative set as it is also based on the mapping on the dual hexagonal lattice.
\par
Finally, all these different sets of objects could be used to generate homochiral inflation rules on their own, a good way to confirm that fixing the chirality of the aperiodic monotile Tile(1,1) is a good simplification.

\section*{Acknowledgements}
We are very grateful for many helpful discussions to Anuradha Jagannathan and Pavel Kalugin at the LPS, Arnaud Chéritat at Institut Mathematiques de Toulouse and Pierre Gradit. All figures representing tiles/clusters are done using the open source software Inkscape (https://inkscape.org). For inflation matrices, Wolfram Mathematica language is used (https://www.wolfram.com/mathematica). 

\section*{References}
1. Smith, D., Myers, J. S., Kaplan, C. S. and Goodman-Strauss, C. A chiral
aperiodic monotile 2023. arXiv: 2305.17743 [math.CO].
\par
2. Baake, M., Gahler, F. and Grimm, U. Hexagonal Inflation Tilings and Planar
Monotiles. Symmetry 4, 581–602. issn: 2073-8994. https://www.mdpi.
com/2073-8994/4/4/581 (2012).
\par
3. Smith, D., Myers, J. S., Kaplan, C. S. and Goodman-Strauss, C. An aperiodic
monotile. arXiv preprint arXiv:2303.10798 (2023).
\par
4. Voss, H. U. A tiling algorithm for the aperiodic monotile Tile(1,1) 2024.
arXiv: 2406.05236 [math-ph]. https://arxiv.org/abs/2406.05236.
\par
5. Akiyama, S. and Araki, Y. An alternative proof for an aperiodic monotile
2025. arXiv: 2307.12322 [math.MG]. https://arxiv.org/abs/2307.
12322.
\par
6. Smith, J. Turtles, Hats and Spectres: Aperiodic structures on a Rhombic
tiling 2024. arXiv: 2403.01911 [math.MG]. https://arxiv.org/abs/
2403.01911.
\par
7. Cheritat, A. Observations on the hex clusters of the Spectre tilings 2024.
arXiv: 2407.05359 [math.CO]. https://arxiv.org/abs/2407.05359.
\par
8. Imperor-Clerc, M., Jagannathan, A., Kalugin, P. and Sadoc, J.-F. Squaretriangle
tilings: an infinite playground for soft matter. Soft Matter 17. Publisher:
The Royal Society of Chemistry, 9560–9575. issn: 1744-6848. https:
//pubs.rsc.org/en/content/articlelanding/2021/sm/d1sm01242h
(2023) (Nov. 3, 2021).
\par
9. Imperor-Clerc, M., Kalugin, P., Schenk, S., Widdra, W. and Forster, S. Higher dimensional
geometrical approach for the characterization of two-dimensional
square-triangle-rhombus tilings. Phys. Rev. B 110, 144106. https://link.
aps.org/doi/10.1103/PhysRevB.110.144106 (14 Oct. 2024).

\begin{figure*}
\includegraphics[width=\textwidth]{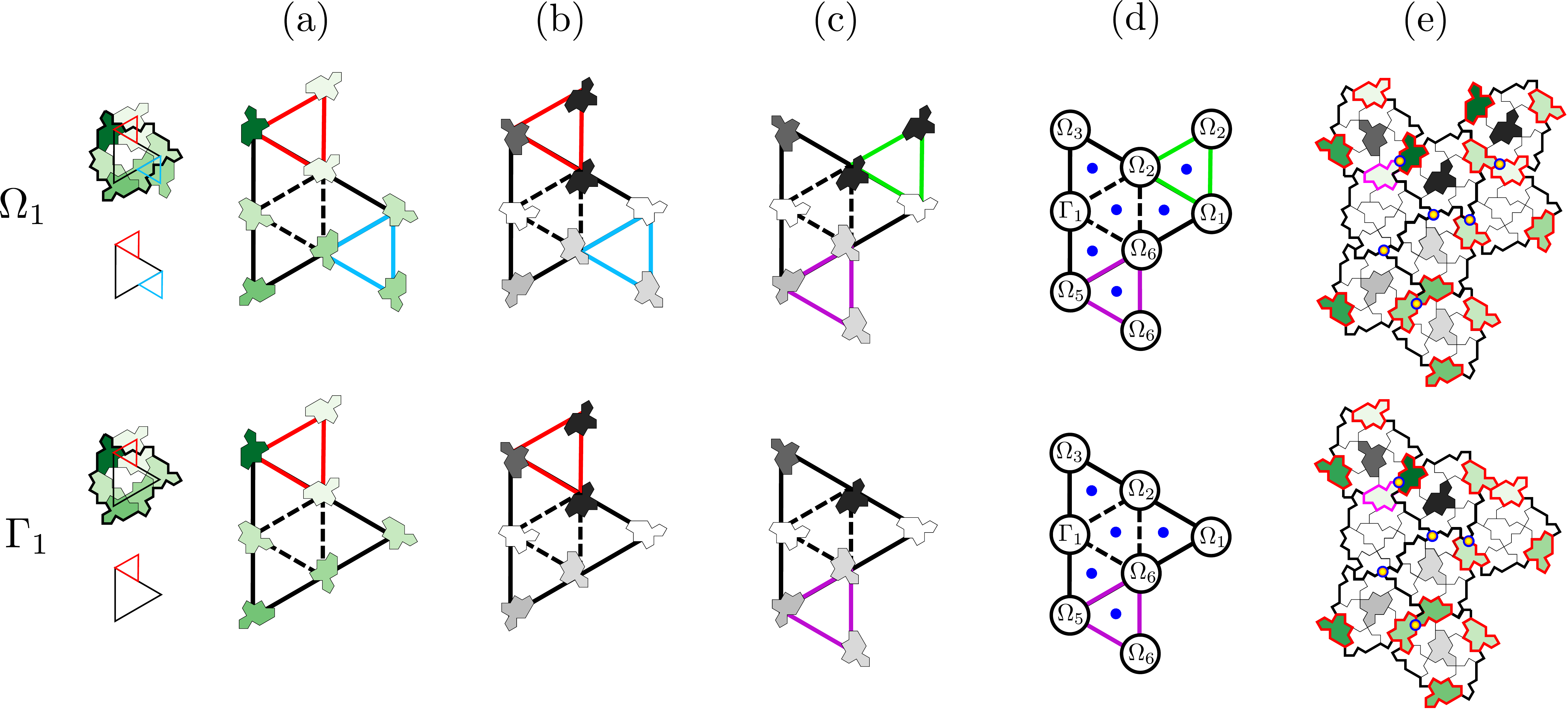}
\caption{Substitution rule for $\Omega_1$ and $\Gamma_1$ metatiles. A triangular pattern is identified in the $\Omega_1$ metatile. The eight even monotiles (shades of green) are first placed on an expanded version of this triangular pattern (a). Then the orientation of these eight monotiles is changed from even to odd by a rotation of $\pi/6$ for all of them (b). Eight odd monotiles (shades of gray) are obtained instead. Finally, the chirality of the triangular pattern is modified as indicated by the change in colors in the triangular patterns (c). This last step involves moving around two positions (red to green triangle and blue to violet triangle in the pattern). The substitution rule of the $\Omega_1$ metatile by 8 metatiles is mapped on a triangular lattice (d) where metatiles are depicted by circles. The corresponding assembly of metatiles (e) is shown. Dark blue dots stand for junction points between three metatiles and superimposed on yellow glue points (see also Figures \ref{fig:translation} and \ref{fig:metatiles}). The substitution rule for the $\Gamma_1$ metatile into 7 metatiles is derived in exactly the same way (bottom at the figure), except that one metatile is taken out from the triangular pattern (green triangle).}
  \label{fig:omega_gamma}
\end{figure*}

\end{document}